

\def\IfDraTex{\if01}

\input amstex
\loadbold
\documentstyle{amsppt}
\IfDraTex
\input DraTex.sty
\input AlDraTex.sty
\def\DraPenSize{0.75pt}
\else
\message{Please use DraTex.sty and AlDraTex.sty.
For details, see the top of this file.}
\fi
\NoBlackBoxes

\pagewidth{32pc}
\pageheight{44pc}
\magnification=\magstep1

\def\QQ{\Bbb{Q}}
\def\RR{\Bbb{R}}
\def\PP{\Bbb{P}}
\def\OO{\Cal{O}}
\def\Supp{\operatorname{Supp}}
\def\Coker{\operatorname{Coker}}
\def\rank{\operatorname{rk}}
\def\Tr{\operatorname{Tr}}
\def\id{\operatorname{id}}
\def\Bs{\operatorname{Bs}}
\def\Hom{\operatorname{Hom}}
\def\Sec{\operatorname{Sec}}
\def\ch{\operatorname{char}}
\def\transdeg{\operatorname{trans.deg}}
\def\rest#1#2{\left.{#1}\right\vert_{{#2}}}

\topmatter
\title
Geometric height inequality on varieties \\
with ample cotangent bundles
\endtitle
\rightheadtext{}
\author Atsushi Moriwaki \endauthor
\leftheadtext{}
\address
Department of Mathematics, Faculty of Science,
Kyoto University, Kyoto, 606-01, Japan
\endaddress
\curraddr
Department of Mathematics, University of California,
Los Angeles, 405 Hilgard Avenue, Los Angeles, California 90024, USA
\endcurraddr
\email moriwaki\@math.ucla.edu \endemail
\date November, 1993 (revised) \enddate
\abstract
Let $F$ be a function field of one variable over
an algebraically closed field of characteristic zero,
$X$ a geometrically irreducible smooth projective variety over $F$,
and $L$ a line bundle on $X$.
In this note, we will prove that if $\Omega^1_{X/F}$ is ample and
$X$ is non-isotrivial, then there are a proper closed
algebraic set $Y$ of $X$ and a constant $A >0$ such that
$$
     h_L(P) \leq A \cdot d(P) + O(1)
$$
for all $P \in X(\bar{F}) \setminus Y(\bar{F})$,
where $h_L(P)$ is a geometric height of $P$
with respect to $L$ and $d(P)$ is
the geometric logarithmic discriminant of $P$.
As corollary of the above height inequality, we can recover
Noguchi's theorem \cite{No}, i.e.
there is a non-empty Zariski open set $U$ of $X$ with $U(F) = \emptyset$.
\endabstract
\endtopmatter

\document

\subhead 0. Introduction
\endsubhead

Let $k$ be an algebraically closed field of characteristic zero
and $F$ a finitely generated extension field of $k$ with $\transdeg_k F = 1$.
We will fix two fields $k$ and $F$ throughout this note.

Let $X$ be a geometrically irreducible smooth projective variety
over $F$.
Let $\Cal{X}$ and $C$ be smooth projective varieties over $k$, and
$f : \Cal{X} \to C$ a $k$-morphism such that the function field
of $C$ is $F$ and the generic fiber of $f$ is $X$, i.e.
$X = \Cal{X} \otimes F$. In this note,
$X$ is said to be {\it non-isotrivial} if
there is a non-empty open set $C_0$ of $C$ such that,
for all $t \in C_0$,
the Kodaira-Spencer map
$\varrho_t : T_{C, t} \to H^1(\Cal{X}_t, T_{\Cal{X}_t})$
is not zero.
Let $\bar{F}$ be the algebraic closure of $F$.
For a point $P \in X(\bar{F})$,
let us denote by $\Delta_P$ the corresponding integral curve on $\Cal{X}$.
We fix a line bundle $L$ on $X$.
Let $\Cal{L}$ be a line bundle on $\Cal{X}$
with $\Cal{L} \otimes F = L$.
The pair $(f : \Cal{X} \to C, \Cal{L})$ is called {\it a model of $(X, L)$}.
{\it A geometric height $h_L(P)$ of $P$ with respect to $L$} is defined by
$$
h_L(P) = \frac{(\Cal{L} \cdot \Delta_P)}{[F(P) : F]}.
$$
The geometric height $h_L$ might depend on the choice of the model
$(f: \Cal{X} \to C, \Cal{L})$ of $(X, L)$, but
as we will see in Proposition~1.1, it is uniquely determined modulo
bounded functions.
Moreover, {\it the geometric logarithmic discriminant $d(P)$ of $P$}
is given by
$$
d(P) = \frac{2 g(\widetilde{\Delta}_P) - 2}{[F(P) : F]},
$$
where $g(\widetilde{\Delta}_P)$ is the genus of the normalization
$\widetilde{\Delta}_P$ of $\Delta_P$.
An interesting problem concerning $h_L$ and $d$ on
the non-isotrivial $X$ is to find an inequality
of type:
$$
    h_L(P) \leq A \cdot d(P) + O(1)
$$
for all $P \in X(\bar{F})$, which is called
{\it a geometric height inequality}.
In the case where $\dim X = 1$,
we know a lot of height inequalities as follows, where $g = g(X) \geq 2$
and $K = K_{\Cal{X}/C}$.
$$\vbox{
\settabs
\+ H. Esnault and E. Viehweg \quad &
$h_{K}(P) \leq (20g - 15)/6 \ d(P) + O(1)$\quad &
($\ch \geq 0$) \cr
\+ L. Szpiro, & $h_{K}(P) \leq 24(g-1)^2 \ d(P) + O(1)$, & ($\ch \geq 0$)
\cr\smallskip
\+ P. Vojta,  & $h_{K}(P) \leq (8g-6)/3 \ d(P) + O(1)$,  & ($\ch =0$)
\cr\smallskip
\+ A. N. Parshin, & $h_{K}(P) \leq (20g - 15)/6 \ d(P) + O(1)$, & ($\ch =0$)
\cr\smallskip
\+ H. Esnault and E. Viehweg, & $h_{K}(P) \leq 2(2g-1)^2 \ d(P) + O(1)$, &
($\ch =0$) \cr\smallskip
\+ P. Vojta, & $h_{K}(P) \leq (2 + \epsilon) \ d(P) + O(1)$, & ($\ch =0$)
\cr\smallskip
\+ A. Moriwaki, & $h_{K}(P) \leq (2g - 1) \ d(P) + O(1)$, & ($\ch \geq 0$)
\cr
}$$
Details can be found in Chapter~VI of \cite{La} or \cite{Mo}.
In this note, inspired by \cite{Vo},
I would like to study a geometric height inequality
on a higher dimensional variety.
The main theorem of this note is

\proclaim{Theorem A}
Let $X$ be a geometrically irreducible smooth projective variety over $F$
and $L$ a line bundle on $X$.
If $\Omega^1_{X/F}$ is ample and
$X$ is non-isotrivial, then there are a proper closed
algebraic set $Y$ of $X$ and a constant $A >0$ such that
$$
     h_L(P) \leq A \cdot d(P) + O(1)
$$
for all $P \in X(\bar{F}) \setminus Y(\bar{F})$.
\endproclaim

Under the assumption of Theorem~A, it might be not needed to
subtract the proper closed algebraic set $Y$.
But, it is sufficient to recover Noguchi's theorem \cite{No}
concerning higher dimensional geometric Mordell's conjecture, i.e.

\proclaim{Theorem B}
Let $X$ be a geometrically irreducible smooth projective variety over $F$.
If $\Omega^1_{X/F}$ is ample and $X$ is non-isotrivial,
then there is a non-empty Zariski open set $U$ of $X$ with
$U(F) = \emptyset$.
\endproclaim

Theorem~B is a corollary of Theorem~A once we get finiteness property
of geometric heights. However, finiteness does not hold in general
for geometric heights. Instead of finiteness, algebraic heights
have non-denseness property (cf. Theorem~1.3).

In \S1, we will deal with properties of geometric heights,
especially `positivity' and `non-denseness'.
In \S2, we will prepare several basic facts of
ample vector bundles.
\S3 is devoted to the proofs of Theorem~A and Theorem~B.

\subhead 1. Properties of geometric heights
\endsubhead

Let $X$ be a geometrically irreducible smooth projective variety over $F$
and $L$ a line bundle on $X$.
First of all, we must check the following proposition.
By this proposition,
$$
h_L(P) = \frac{(\Cal{L} \cdot \Delta_P)}{[F(P) : F]}
$$
is uniquely determined modulo bounded functions.

\proclaim{Proposition 1.1}
Let $(f : \Cal{X} \to C, \Cal{L})$ and $(f' : \Cal{X}' \to C, \Cal{L}')$
be two models of $(X, L)$. Then, there is a constant $A$ such that,
for all $P \in X(\bar{F})$,
$$
 \left| \frac{(\Cal{L} \cdot \Delta_P)}{[F(P) : F]} -
        \frac{(\Cal{L}' \cdot \Delta'_P)}{[F(P) : F]} \right|
 \leq A.
$$
\endproclaim

\demo{Proof}
We can easily construct birational morphisms
$\mu : \Cal{X}'' \to \Cal{X}$ and $\mu' : \Cal{X}'' \to \Cal{X}'$
such that $(\Cal{X}'' \to C, \mu^*(\Cal{L}))$ and
$(\Cal{X}'' \to C, {\mu'}^*(\Cal{L}'))$ are models of $(X, L)$.
So we may assume that $\Cal{X} = \Cal{X}'$.

Since $\Cal{L} \otimes {\Cal{L}'}^{-1}$ is trivial on $X$,
we have $f_*(\Cal{L} \otimes {\Cal{L}'}^{-1}) \not= 0$.
Thus there is an ample line bundle $A_1$ on $C$ such that
$H^0(\Cal{X}, \Cal{L} \otimes {\Cal{L}'}^{-1} \otimes f^*(A_1)) \not= 0$.
Let $D$ be an element of
$|\Cal{L} \otimes {\Cal{L}'}^{-1} \otimes f^*(A_1)|$.
$D$ has no horizontal components because
$\rest{D}{X}$ is linearly equivalent to zero. Thus,
for all $P \in X(\bar{F})$, $(D \cdot \Delta_P) \geq 0$.
It follows that
$$
\frac{(\Cal{L}' \cdot \Delta_P)}{[F(P) : F]} -
        \frac{(\Cal{L} \cdot \Delta_P)}{[F(P) : F]} \leq \deg(A_1).
$$
By the same way, we can find an ample line bundle $A_2$ on $C$ such that
$$
\frac{(\Cal{L} \cdot \Delta_P)}{[F(P) : F]} -
        \frac{(\Cal{L}' \cdot \Delta_P)}{[F(P) : F]} \leq \deg(A_2).
$$
Hence, we get our proposition.
\qed
\enddemo

Let $V$ be a projective variety defined over an algebraic number
field $K$, $L$ a line bundle on $V$, and $h_L$ a Weil-height of $L$.
Then, $h_L$ has the following two basic properties:
\roster
\item "(1)" (Positivity) Let $B$ be a base locus of $|L|$.
Then, $h_L(P) \geq O(1)$ for all
$P \in V(\overline{\QQ}) \setminus B(\overline{\QQ})$.

\item "(2)" (Finiteness) If $L$ is ample and $A \in \RR$, then
the set $\{ P \in V(K) \mid h_L(P) \leq A \}$ is finite.
\endroster
The main topic of this section is to consider geometric analogies
of the above two properties of arithmetic heights.
First, let us consider positivity of geometric heights.

\proclaim{Proposition 1.2}(Positivity)\quad
Let $X$ be a geometrically irreducible smooth projective variety over $F$
and $L$ a line bundle on $X$.
Then, $h_L(P) \geq O(1)$
for all $P \in X(\bar{F}) \setminus \Bs(L)(\bar{F})$, where
$\Bs(L) = \Supp(\Coker(H^0(X, L) \otimes \OO_X \to L))$.
\endproclaim

\demo{Proof}
Let $(f : \Cal{X} \to C, \Cal{L})$ be a model of $(X, L)$.
Let $A$ be an ample line bundle on $C$ such that
$f_*(\Cal{L}) \otimes A$ is generated by global sections.
Let $P \in X(\bar{F}) \setminus \Bs(L)(\bar{F})$ and
$\Delta_P$ the corresponding integral curve on $\Cal{X}$.
Assume that $\Delta_P \subset \Bs(\Cal{L} \otimes f^*(A))$.
Then, the natural homomorphism
$$
    H^0(C, f_*(\Cal{L} \otimes f^*(A) \otimes I_{\Delta_P})) \to
    H^0(C, f_*(\Cal{L} \otimes f^*(A)))
$$
is surjective, where $I_{\Delta_P}$ is the defining ideal of $\Delta_P$.
Here, we consider the following commutative diagram:
$$
\CD
    H^0(C, f_*(\Cal{L} \otimes f^*(A) \otimes I_{\Delta_P}))\otimes \OO_C
    @>{\alpha}>>
    H^0(C, f_*(\Cal{L} \otimes f^*(A))) \otimes \OO_C \\
    @V{\gamma}VV   @VV{\delta}V \\
    f_*(\Cal{L} \otimes f^*(A) \otimes I_{\Delta_P})
    @>>{\beta}>
    f_*(\Cal{L} \otimes f^*(A))
\endCD
$$
Since $\alpha$ and $\delta$ are surjective,
so is $\beta$, which implies that
$H^0(Y_{\bar{F}}, L_{\bar{F}} \otimes I_P) \to H^0(Y_{\bar{F}}, L_{\bar{F}})$
is surjective. Therefore, $P \in \Bs(L_{\bar{F}})$.
This is a contradiction. Thus, we get
$\Delta_P \not\subset \Bs(\Cal{L} \otimes f^*(A))$.
Hence, $(\Cal{L} \otimes f^*(A) \cdot \Delta_P) \geq 0$, which says that
$h_L(P) \geq -\deg A$.
\qed
\enddemo

Next, we will consider finiteness of geometric heights.
Unfortunately, finiteness does not hold in general for geometric heights.
Instead of it, we have the following non-denseness of geometric heights.

\proclaim{Theorem 1.3}(Non-denseness)\quad
Let $X$ be a geometrically irreducible smooth projective variety over $F$
and $L$ an ample line bundle on $X$.
We assume that $X$ contains no rational curves and
$X$ is non-isotrivial.
Then, for every constant $A$,
the set $\{ P \in X(F) \mid h_L(P) \leq A \}$
is not dense in $X(\bar{F})$.
\endproclaim

First of all, we will prepare two basic lemmas.

\proclaim{Lemma 1.4}
Let $q : Y \to X$ be a surjective morphism of normal projective
varieties over an algebraically closed field
of characteristic zero and $E$ a vector bundle on $X$. Then, the natural
morphism $H^1(X, E) \to H^1(Y, q^* E)$ is injective.
\endproclaim

\demo{Proof}
Let us consider the Leray spectral sequence:
$$
    E^{i,j}_2 = H^i(X, R^iq_*(q^* E)) \Longrightarrow
    E^{i+j} = H^{i+j}(Y, q^* E).
$$
Since
$$
    E^{1,0}_2 = H^1(X, q_*(q^* E)) \to
    E^1 = H^1(Y, q^* E)
$$
is injective, it is sufficient to see that
$$
    H^1(X, E) \to H^1(X, q_*(\OO_Y) \otimes E)
$$
is injective. Let
\IfDraTex
\par\medskip
\Draw
\PenSize(\DraPenSize)
\Node(a)(--$Y$--)
\Move(70,0)
\Node(b)(--$Z$--)
\Move(70,0)
\Node(c)(--$X$--)
\ArrowHeads(1) \ArrowSpec(V)
\Edge(a,b) \EdgeLabel[+](--$q_1$--)
\Edge(b,c) \EdgeLabel[+](--$q_2$--)
\CurvedEdgeSpec(20,0.2,-20,0.2)
\CurvedEdge(a,c) \EdgeLabel(--$q$--)
\EndDraw
\medskip\par\noindent
\else
$$
\CD
 Y @>{q_1}>> Z @>{q_2}>> X
\endCD
$$
\fi
be the Stein factorization of $q : Y \to X$.
Then, $q_*(\OO_Y) = {q_2}_*(\OO_Z)$. Thus, we may assume that
$q$ is a finite morphism.
Since we work over the field of characteristic zero,
we have a trace map
$$
    \Tr : q_* \OO_Y \to \OO_X
$$
such that $\Tr \cdot \iota = \id$, where $\iota$ is the inclusion map
$\iota : \OO_X \to q_*(\OO_Y)$. This means that
the exact sequence
$$
   0 \to E \to q_*(\OO_Y) \otimes E \to
  \left( q_*(\OO_Y)/\OO_X \right) \otimes E \to 0
$$
splits. Therefore,
$H^1(X, E) \to H^1(X, q_*(\OO_Y) \otimes E)$ is injective.
\qed
\enddemo

\proclaim{Lemma 1.5}
Let $\Cal{X}$ be a smooth projective variety over $k$,
$C$ a smooth projective curve over $k$,
$f : \Cal{X} \to C$ a $k$-morphism with $f_* \OO_\Cal{X} = \OO_C$.
Assume that there are a non-empty open set $C_0$ of $C$,
a smooth projective variety $Y$ over $k$ and a dominant morphism
$\phi : Y \times C_0 \to \Cal{X}_0 = f^{-1}(C_0)$ such that
$p = f_0 \cdot \phi$, where $p : Y \times C_0 \to C_0$ is the natural
projection and $f_0 = \rest{f}{\Cal{X}_0}$.
\IfDraTex
\par\medskip
\Draw
\PenSize(\DraPenSize)
\Node(a)(--$Y \times C_0$--)
\Move(70,0)
\Node(b)(--$\Cal{X}_0$--)
\Move(-35,-50)
\Node(c)(--$C_0$--)
\ArrowHeads(1) \ArrowSpec(V)
\Edge(a,b) \EdgeLabel(--$\phi$--)
\Edge(b,c) \EdgeLabel(--$f_0$--)
\Edge(a,c) \EdgeLabel[+](--$p$--)
\EndDraw
\medskip\par\noindent
\else
$$
\CD
        Y \times C_0 @>{\phi}>> \Cal{X}_0        \\
        @V{p}VV                 @VV{f_0}V  \\
        C_0          @=         C_0
\endCD
$$
\fi
Then, for a general point $t$ of $C_0$,
the Kodaira-Spencer map
$\varrho_t : T_{C, t} \to H^1(\Cal{X}_t, T_{\Cal{X}_t})$ is zero.
\endproclaim

\demo{Proof}
Shrinking $C_0$, if necessarily, we may assume that $f_0$ is smooth.
Let us consider the following commutative diagram:
$$
\CD
0 @>>> T_{Y \times C_0/C_0} @>>>T_{Y \times C_0} @>>>
       p^*(T_{C_0}) @>>> 0 \\
@. @VVV @VVV @| @. \\
0 @>>> \phi^*(T_{\Cal{X}_0/C_0}) @>>> \phi^*(T_{\Cal{X}_0}) @>>>
       \phi^*(f_0^*(T_{C_0})) @>>> 0
\endCD
$$
This gives rise to a commutative diagram:
\IfDraTex
\par\medskip
\Draw
\PenSize(\DraPenSize)
\Node(a)(--$T_{C,t}$--)
\Move(100,30)
\Node(b)(--$H^1(Y \times \{ t \}, T_{Y \times \{ t \}})$--)
\Move(0,-60)
\Node(c)(--$H^1(Y \times \{ t \}, \phi_t^*(T_{\Cal{X}_t}))$--)
\ArrowHeads(1) \ArrowSpec(V)
\Edge(a,b) \EdgeLabel(--$\alpha_t$--)
\Edge(b,c)
\Edge(a,c) \EdgeLabel[+](--$\beta_t$--)
\EndDraw
\medskip\par\noindent
\else
$$
\CD
T_{C, t} @>{\alpha_t}>> H^1(Y \times \{ t \}, T_{Y \times \{ t \}}) \\
@|                        @VVV  \\
T_{C, t} @>>{\beta_t}>  H^1(Y \times \{ t \}, \phi_t^*(T_{\Cal{X}_t}))
\endCD
$$
\fi
for a point $t \in C_0$.
Since $\alpha_t = 0$, we have $\beta_t = 0$.
Here $\beta_t$ is factored through the Kodaira-Spencer map
$\varrho_t : T_{C, t} \to H^1(\Cal{X}_t, T_{\Cal{X}_t})$, i.e.
if we denote by $\gamma_t$ the natural homomorphism
$ H^1(\Cal{X}_t,T_{\Cal{X}_t}) \to
H^1(Y \times \{ t \}, \phi_t^*(T_{\Cal{X}_t}))$, then
$\beta_t = \gamma_t \cdot \varrho_t$.
\IfDraTex
\par\medskip
\Draw
\PenSize(\DraPenSize)
\Node(a)(--$T_{C,t}$--)
\Move(75,0)
\Node(b)(--$H^1(\Cal{X}_t, T_{\Cal{X}_t})$--)
\Move(115,0)
\Node(c)(--$H^1(Y \times \{ t \}, \phi_t^*(T_{\Cal{X}_t}))$--)
\ArrowHeads(1) \ArrowSpec(V)
\Edge(a,b) \EdgeLabel[+](--$\varrho_t$--)
\Edge(b,c) \EdgeLabel[+](--$\gamma_t$--)
\CurvedEdgeSpec(20,0.2,-20,0.2)
\CurvedEdge(a,c) \EdgeLabel(--$\beta_t$--)
\EndDraw
\medskip\par\noindent
\else
$$
\CD
 T_{C, t} @>{\varrho_t}>> H^1(\Cal{X}_t, T_{\Cal{X}_t}) @>{\gamma_t}>>
 H^1(Y \times \{ t \}, \phi_t^*(T_{\Cal{X}_t}))
\endCD
$$
\fi
On the other hand, by Lemma~1.4,
$\gamma_t$ is injective. Thus, $\varrho_t = 0$
\qed
\enddemo

\demo{Proof of Theorem~1.3}
Let $(f : \Cal{X} \to C, \Cal{L})$ be a model of $(X, L)$.
Let $\Cal{M}$ be an ample line bundle on $\Cal{X}$.
Since $L$ is ample on $X$, for a sufficiently large $n$,
$\Cal{L}^n \otimes \Cal{M}^{-1}$ is very ample on $X$.
Thus, by Proposition~1.1, there is a constant $B$ such that
$h_{L^n \otimes M^{-1}}(P) \geq B$ for all $P \in X(\bar{F})$,
where $M = \Cal{M}_F$.
Therefore, we have
$$
     h_M(P) \leq n \cdot h_L(P) - B
$$
for all $P \in X(\bar{F})$. It follows that we may assume that
$\Cal{L}$ is ample on $\Cal{X}$.

We set $S_A  = \{ P \in X(F) \mid  h_L(P) \leq A \}$ and
$\Sigma_A = \{ \Delta_P \subset \Cal{X} \mid P \in S_A \}$.
Here we assume that
$S_A$ is dense in $X(\bar{F})$.
First of all, we claim that
there is a smooth quasi-projective variety $Y$ over $k$ and a dominant morphism
$\phi : Y \times C \to \Cal{X}$ such that
$p = f \cdot \phi$, where $p : Y \times C \to C$ is the natural
projection.
\IfDraTex
\par\medskip
\Draw
\PenSize(\DraPenSize)
\Node(a)(--$Y \times C$--)
\Move(70,0)
\Node(b)(--$\Cal{X}$--)
\Move(-35,-50)
\Node(c)(--$C$--)
\ArrowHeads(1) \ArrowSpec(V)
\Edge(a,b) \EdgeLabel(--$\phi$--)
\Edge(b,c) \EdgeLabel(--$f$--)
\Edge(a,c) \EdgeLabel[+](--$p$--)
\EndDraw
\medskip
\else
$$
\CD
        Y \times C @>{\phi}>> \Cal{X}        \\
        @V{p}VV               @VV{f}V  \\
        C          @=         C
\endCD
$$
\fi

Let $\Hom(C, \Cal{X})$ be a scheme consisting of $k$-morphisms from $C$ to
$\Cal{X}$. Then, we have the natural morphism of schemes
$\alpha : \Hom(C, \Cal{X}) \to \Hom(C, C)$
defined by $\alpha(s) = f \cdot s$
for $s \in \Hom(C, \Cal{X})$. We set
$\Sec(C, \Cal{X}) = \alpha^{-1}(\id_C)$.

By the definition of $\Sigma_A$, $\Sigma_A$ is a bounded family.
So we can find finite number of
connected components $\Sec_1, \ldots, \Sec_l$ of $\Sec(C, \Cal{X})$
such that, for all $P \in S_A$,
there is $s \in \coprod_{i=1}^{l} \Sec_i$ with $s(C) = \Delta_P$.
Further, since every connected component of $\Hom(C, \Cal{X})$
is quasi-projective, $\Sec_i$'s are also quasi-projective.

On the other hand, since $S_A$ is dense in $X(\bar{F})$,
there is a $\Sec_i$ such that the natural morphism
$\Sec_i \times C \to \Cal{X}$
is a dominant morphism. Moreover, since $\Cal{X}$ is irreducible,
there is an irreducible component $Y$ of $\Sec_i$ such that
$Y \times C \to \Cal{X}$ still dominates $\Cal{X}$.
Replacing $Y$ by $Y_{\operatorname{red}}$, we may assume that $Y$
is reduced. Furthermore, considering a desingularization $Y' \to Y$
of $Y$, we may also assume that $Y$ is smooth.
Thus, we have our claim.

Let $Y_1$ be a smooth compactification of $Y$.
Then, $\phi$ induces a rational map
$\phi_1 : Y_1 \times C \dashrightarrow \Cal{X}$.
Let $\mu : Z \to Y_1 \times C$ be a minimal succession of blowing-ups such that
$\phi_1 \cdot \mu : Z \to \Cal{X}$ gives a morphism.
Let $E$ be the exceptional set of $\mu$.
Since $X$ contains no rational curves,
$(p \cdot \mu)(E)$ is a proper closed subset of $C$.
Therefore, we can take a non-empty open set $C_0$
of $C$ such that $\phi_1$ is defined on $Y_1 \times C_0$.
Thus, by virtue of Lemma~1.5,
for a general point $t$ of $C_0$,
the Kodaira-Spencer map
$\varrho_t : T_{C, t} \to H^1(\Cal{X}_t, T_{\Cal{X}_t})$ is zero.
This is a contradiction.
Therefore, $S_A$ is not dense in $X(\bar{F})$.
\qed
\enddemo

\definition{Remark 1.6}
For geometric heights, finiteness does not hold in general.
For, as pointing out in
Example of \cite{No}, we have a non-isotrivial family of
projective varieties with ample cotangent bundles, which contains
a trivial family of subvarieties. More precisely,
we have the following example with properties (1) --- (5).

\roster
\item "(1)" Let $f : \Cal{X} \to C$ be a projective morphism
of smooth varieties.

\item "(2)" $\dim \Cal{X} = 3$, $\dim C = 1$ and $f_*\OO_{\Cal{X}} = \OO_C$.

\item "(3)" $\Omega^1_{\Cal{X}_t}$ is ample for a general $t \in C$.

\item "(4)" The Kodaira-Spencer map
$\varrho_t : T_{C, t} \to H^1(\Cal{X}_t, T_{\Cal{X}_t})$
is injective for a general $t \in C$.

\item "(5)" There are a smooth projective curve $C'$ and
an embedding $j : C' \times C \to \Cal{X}$ with $f \cdot j = p$,
where $p$ is the natural projection $p : C' \times C \to C$.
\endroster
\IfDraTex
\par\medskip
\Draw
\PenSize(\DraPenSize)
\Node(a)(--$C' \times C$--)
\Move(70,0)
\Node(b)(--$\Cal{X}$--)
\Move(-35,-50)
\Node(c)(--$C$--)
\ArrowHeads(1) \ArrowSpec(V)
\Edge(a,b) \EdgeLabel(--$j$--)
\Edge(b,c) \EdgeLabel(--$f$--)
\Edge(a,c) \EdgeLabel[+](--$p$--)
\EndDraw
\medskip
\else
$$
\CD
C' \times C @>{j}>> \Cal{X} \\
@V{p}VV             @VV{f}V \\
C           @=      C
\endCD
$$
\fi
\enddefinition

\subhead 2. Ample vector bundles
\endsubhead

In this section, we will prepare basic facts of ample
vector bundles.
Let us begin with Hartshorne's criterion of ample vector bundles on
curves.

\proclaim{Lemma 2.1}(cf. \cite{Ha, Theorem 2.4})\quad
Let $C$ be a smooth projective curve over
an algebraically closed field of characteristic zero and $E$
a vector bundle on $C$ with $\deg(E) > 0$.
Then, $E$ is ample if and only if
every quotient bundle of $E$ has positive degree.
\endproclaim

\demo{Proof}
For reader's convenience, we will give a sketch of an elementary proof.
Clearly, it is sufficient to show that
if every quotient bundle of $E$ has positive degree,
then $E$ is ample.

First, let us consider a case where $E$ is semistable.
Let $\pi : \PP(E) \to C$ be the projective bundle of $E$ and
$\OO(1)$ the tautological line bundle of $E$.
Then, by \cite{Mi, Theorem~3.1},
$\OO(r) \otimes \pi^*(\det(E)^{-1})$ is nef,
where $r = \rank E$.
Thus, $\OO(r)$ is ample because $\deg (E) > 0$.

Next, we deal with a general case.
We prove it by induction on $\rank E$.
By the previous observation, we may assume that $E$ is not semistable.
Let $E_1$ be the maximal destabilizing subbundle of $E$.
Since $E_1$ is semistable and $\deg (E_1) > 0$, $E_1$ is ample
as before. On the other hand, by hypothesis of induction, $E/E_1$ is ample.
Thus, $E$ is ample.
\qed
\enddemo

\proclaim{Corollary 2.2}
Let $C$ be a smooth projective curve over
an algebraically closed field of characteristic zero and $R$
an ample vector bundle on $C$.
Let
$$
0 \to \OO_C \to E \to R \to 0
\tag 2.2.1
$$
be an extension of $R$ by $\OO_C$.
Then, $E$ is ample if and only if
(2.2.1) does not split.
\endproclaim

\demo{Proof}
First, we assume that $E$ is ample.
Let $\pi : \PP(E) \to C$ be the projective bundle of $E$
and $\OO(1)$ the tautological line bundle.
If (2.2.1) splits, then we have a section $\Delta$ of $\pi$ defined by
a splitting $E \to \OO_C$. Clearly, $(\OO(1) \cdot \Delta) = 0$.
This is a contradiction because $\OO(1)$ is ample.
Hence, (2.2.1) doesn't split.

Next we assume that (2.2.1) doesn't split.
Let $Q$ be a quotient bundle of $E$.
Let $\alpha : \OO_C \to E \to Q$ be the composition homomorphism.
If $\alpha$ is zero, we have a surjective homomorphism $R \to Q$.
Thus, $\deg(Q) > 0$. So we may assume that $\alpha \not= 0$.
Let $L$ be the saturation of $\alpha(\OO_C)$ in $Q$.
Let us consider the following homomorphisms:
$$
R \simeq E/\OO_C \to Q/\alpha(\OO_C) \to (Q/\alpha(\OO_C))/\hbox{(torsions)}
\simeq Q/L.
$$
This shows that there is a surjective homomorphism $R \to Q/L$.
Thus, if $\rank Q \geq 2$, then
$\deg(Q) = \deg(Q/L) + \deg(L) > 0$.
Moreover, in the case where $\rank Q = 1$,
since (2.2.1) doesn't split, $\alpha$ is not an
isomorphism. Hence, $\deg(Q) = \deg(L) > 0$.
Therefore, $E$ is ample by Lemma~2.1.
\qed
\enddemo

\proclaim{Proposition 2.3}
Let $X$ be a smooth projective variety over an algebraically closed field
of characteristic zero and
$Q$ an ample vector bundle on $X$.
Let $0 \to \OO_X \to E \to Q \to 0$ be an
extension of $Q$ by $\OO_X$.
Let $\pi : P = \PP(E) \to X$ be the projective bundle of $E$,
$\OO_P(1)$ the tautological line bundle of $P$, and $M$ a line bundle
on $P$. If $0 \to \OO_X \to E \to Q \to 0$ does not
split, then, for a sufficiently large $n$,
$\pi(\Bs(\OO_P(n) \otimes M))$ is a proper closed subset of $X$.
\endproclaim

\demo{Proof}
Let $Y$ be a sub-variety of $P$ defined by $E \to Q$.
Then, $Y \in |\OO_P(1)|$. First of all, since $\rest{\OO_P(1)}{Y}$
is ample, if $n \gg 0$, then
$H^i(Y, \rest{\OO_P(n) \otimes M}{Y}) = 0$ for $i > 0$ and
$\rest{\OO_P(n) \otimes M}{Y}$ is generated by global sections.
Let us consider an exact sequence:
$$
0 \to \OO_P(n-1) \otimes M \to \OO_{P}(n) \otimes M \to
      \rest{\OO_P(n) \otimes M}{Y} \to 0.
$$
If $n \gg 0$,
$$
  H^1(P, \OO_P(n-1) \otimes M) \to H^1(P, \OO_P(n) \otimes M)
\tag 2.3.1
$$
is surjective.
So, $h^1(P, \OO_P(n) \otimes M)$ is constant for $n \gg 0$.
Therefore, (2.3.1) is bijective for $n \gg 0$.
Hence, we get
$$
     H^0(P, \OO_P(n) \otimes M) \to H^0(Y, \rest{\OO_P(n) \otimes M}{Y})
$$
is surjective for $n \gg 0$. It follows that, if we set
$B_n = \Bs(\OO_P(n) \otimes M)$, then $B_n \cap Y = \emptyset$
for $n \gg 0$.

Let $H$ be an ample line bundle on $X$.
If $\dim X \geq 2$, then $H^1(X, Q^{\vee} \otimes H^{-m}) = 0$
for sufficiently large $m$.
Thus, choosing a general member
$T$ of $|H^m|$, the restriction map
$H^1(X, Q^{\vee}) \to H^1(T, \rest{Q^{\vee}}{T})$
is injective. This observation shows us that,
for a general curve $C$ on $X$,
$H^1(X, Q^{\vee}) \to H^1(C, \rest{Q^{\vee}}{C})$
is injective. This means that the exact sequence:
$$
0 \to \OO_C \to \rest{E}{C} \to \rest{Q}{C} \to 0
$$
does not split for the curve $C$.
Hence, by Corollary~2.2, $\rest{E}{C}$ is ample.

Here we assume that $\pi(B_n) = X$.
Then, there is an integral curve $C'$ such that $C' \subset B_n$ and
$\pi(C') = C$. Then, since $C' \cap Y = \emptyset$, we have
$(\OO_P(1) \cdot C') = 0$. This contradicts to ampleness of $\rest{E}{C}$.
\qed
\enddemo

\subhead 3. Proofs of Theorem~A and Theorem~B
\endsubhead

\subsubhead\nofrills{\bf 3.1. Proof of Theorem~A.}
\endsubsubhead\quad
Let $(f : \Cal{X} \to C, \Cal{L})$ be a model of $(X, L)$.
Let $\pi : \Cal{Y} = \PP(\Omega^1_{\Cal{X}/k}) \to \Cal{X}$ be
the projective bundle and
$\OO_\Cal{Y}(1)$ the tautological line bundle on $\Cal{Y}$.
We set $Y = \Cal{Y}_F$, $B'_n = \Bs(\OO_{Y}(n) \otimes \pi_{F}^*(L^{-1}))$,
and $B_n = \pi_{F}(B'_n)$.
By Proposition~1.2, there is a constant $D$ such that, for all
$P' \in Y(\bar{F}) \setminus B'_n(\bar{F})$,
$$
    h_{\OO_{Y}(n) \otimes \pi_F^*(L^{-1})}(P') \geq D.
\tag 3.1.1
$$

On the other hand, the canonical exact sequence:
$$
  0 \to f^*(\Omega^1_{C/k}) \to \Omega^1_{\Cal{X}/k} \to
  \Omega^1_{\Cal{X}/C} \to 0
$$
gives rise to an exact sequence on $X$:
$$
  0 \to \OO_{X} \to \rest{\Omega^1_{\Cal{X}/k}}{X} \to
        \Omega^1_{X/F} \to 0.
\tag 3.1.2
$$
Since $X$ is non-isotrivial, (3.1.2) does not split.
Thus, by virtue of Proposition~2.3, if $n \gg 0$, then $B_n \otimes \bar{F}$
is a proper closed subset of $X_{\bar{F}}$. Therefore, so is
$B_n$.

Let $P \in X(\bar{F}) \setminus B_n(\bar{F})$, $\Delta_P$
the corresponding integral curve on $\Cal{X}$,
$\widetilde{\Delta}_P$ the normalization of
$\Delta_P$ and $\iota : \widetilde{\Delta}_P \to \Cal{X}$ the
natural morphism.
Moreover, let $\phi : \PP(\iota^*(\Omega^1_{\Cal{X}/k})) \to \Cal{Y}$
be the induced morphism.
$$
\CD
   \PP(\iota^*(\Omega^1_{\Cal{X}/k})) @>{\phi}>>
   \Cal{Y} = \PP(\Omega^1_{\Cal{X}/k}) \\
   @VVV                                @VV{\pi}V \\
   \widetilde{\Delta}_P     @>>{\iota}>      \Cal{X}
\endCD
$$
Let $Q$ be the minimal quotient line bundle of
$\iota^*(\Omega^1_{\Cal{X}/k})$ and
$T$ the corresponding section of $\PP(\iota^{*}(\Omega^1_{\Cal{X}/k}))$.
Let $\Delta' = \phi(T)$ and $P'$ the corresponding $\bar{F}$-value
point of $Y$. Then, $P' \in Y(\bar{F}) \setminus B'_n(\bar{F})$.
Thus, by (3.1.1),
$$
    D \leq h_{\OO_{Y}(n) \otimes \pi_F^*(L^{-1})}(P') =
    n \frac{\deg(Q)}{[F(P) : F]} - h_L(P).
\tag 3.1.3
$$
On the other hand, since there is a non-trivial homomorphism
$\iota^*(\Omega^1_{\Cal{X}/k}) \to \Omega^1_{\widetilde{\Delta}_P/k}$,
we have $\deg Q \leq \deg (\Omega^1_{\widetilde{\Delta}_P/k})$.
Thus, by (3.1.3), we get
$h_L(P) \leq n \cdot d(P) - D$.
\qed

\subsubhead\nofrills{\bf 3.2. Proof of Theorem~B.}
\endsubsubhead\quad
By Theorem~A and Theorem~1.3, it is sufficient to see the
following lemma.

\proclaim{Lemma 3.2.1}
Let $X$ be a smooth projective variety over
a field $K$ with the ample cotangent bundle and
$C$ a $1$-dimensional subvariety of $X$.
Then, the genus of $C$ is greater than one.
\endproclaim

\demo{Proof}
Let $\widetilde{C}$ be the normalization of $C$ and
$\phi : \widetilde{C} \to X$ the induced morphism.
Then, we have a non-trivial homomorphism
$\phi^*(\Omega^1_{X/K}) \to \Omega^1_{\widetilde{C}/K}$.
Here, $\phi^*(\Omega^1_{X/K})$ is ample.
Therefore, $\deg(\Omega^1_{\widetilde{C}/K}) > 0$.
\qed
\enddemo



\widestnumber\key{BPV}
\Refs

\ref\key Ha
\by R. Hartshorne
\paper Ample vector bundles on curves
\jour Nagoya Math. J.
\vol 43
\yr 1971
\pages 73--89
\endref

\ref\key La
\by S. Lang (Ed.)
\book Encyclopaedia of Mathematical Sciences, Vol. 60, Number Theory III,
\yr 1991
\publ Springer-Verlag
\endref

\ref\key Mi
\by Y. Miyaoka
\paper The Chern classes and Kodaira dimension of a minimal variety
\jour Advanced Studies in Pure Mathematics 10, 1987,
Algebraic Geometry, Sendai, 1985
\pages 449--476
\publ Kinokuniya
\endref

\ref\key Mo
\by A. Moriwaki
\paper Height inequality of non-isotrivial curves over function fields
\jour to appear in J. Algebraic Geometry
\endref

\ref\key No
\by J. Noguchi
\paper A higher dimensional analogue of Mordell's
conjecture over function fields
\jour Math. Ann.
\vol 258
\yr 1981
\pages 207--212
\endref

\ref\key Vo
\by P. Vojta
\paper On algebraic points on curves
\jour Comp. Math.
\vol 78
\yr 1991
\pages 29--36
\endref

\endRefs
\enddocument